\documentclass[letterpaper, 10 pt, conference]{ieeeconf}

\IEEEoverridecommandlockouts 
\usepackage{amsmath, amssymb, mathtools, mathrsfs, xcolor}
\usepackage{pgfplots}
\usepackage{tikz}

\pgfplotsset{every axis/.append style={
		tick label style={font=\scriptsize}  
}}

\DeclareMathOperator{\im}{im}

\newcommand{\bbR}{\mathbb{R}}

\newcommand{\bbN}{\mathbb{N}}
\newcommand{\bbZ}{\mathbb{Z}}

\newcommand{\bbT}{\mathbb{T}}

\newcommand{\B}{\mathscr{B}}

\newcommand{\calH}{\mathcal{H}}

\newcommand{\calQ}{\mathcal{Q}}

\newcommand{\calR}{\mathcal{R}}

\newcommand{\calU}{\mathcal{U}}
\newcommand{\calW}{\mathcal{W}}

\newcommand{\bbm}{\begin{bmatrix*}}
	\newcommand{\ebm}{\end{bmatrix*}}

\renewcommand{\t}{^\top}

\newcommand{\norm}[1]{\lVert#1\rVert}
\newcommand{\set}[2]{\left\{ #1 \,\left|\, \vphantom{#1} #2 \right. \right\}}

\newcommand{\inner}[2]{\langle#1, #2\rangle}

\newcommand{\oneto}[1]{[#1]}
\newcommand{\half}{\frac{1}{2}}

\newcommand{\hatQ}{\hat\calQ}
\newcommand{\Gtrue}{G_{\textup{true}}}

\newtheorem{example}{Example}
\newtheorem{definition}{Definition}
\newtheorem{lemma}{Lemma}
\newtheorem{theorem}{Theorem}
\newtheorem{proposition}{Proposition}
\newtheorem{problem}{Problem}

\newtheorem{corollary}{Corollary}

\title{\LARGE \bf
	Towards a representer theorem for identification of passive systems
}

\author{B. M. Shali \quad and \quad H. J. van Waarde% <-this % stops a space
\thanks{The authors B. M. Shali (\texttt{b.m.shali@rug.nl}) and H. J. van Waarde (\texttt{h.j.van.waarde@rug.nl}) are with the Bernoulli Institute for Mathematics, Computer Science, and Artificial Intelligence, University of Groningen, Groningen, The Netherlands. H. J. van Waarde acknowledges financial support by the Dutch Research Council under the NWO Talent Programme Veni Agreement (VI.Veni.22.335).
}% <-this % stops a space
}

\begin{document}

\maketitle
\thispagestyle{empty}
\pagestyle{empty}

\begin{abstract}
		A major problem in system identification is the incorporation of prior knowledge about the physical properties of the given system, such as stability, positivity and passivity. In this paper, we present first steps towards tackling this problem for passive systems. In particular, using ideas from the theory of reproducing kernel Hilbert spaces, we solve the problem of identifying a nonnegative input-output operator from data consisting of input-output trajectories of the system. We prove a representer theorem for this problem in the case where the input space is finite-dimensional. This provides a computationally tractable solution, which we show can be obtained by solving an associated semidefinite program.
\end{abstract}

\section{Introduction}

	System identification is the process of obtaining a mathematical model that accurately describes the behaviour of a dynamical system based on empirical observations. System identification has a rich history in the field of systems and control. The classical approach is to search within a given class of finite-dimensional parametric models, e.g., in autoregressive or state-space form. Recently, inspired by developments in machine learning, a new paradigm for system identification has emerged \cite{pilonetto2011, pillonetto2014, dinuzzo2015, ljung2020}. The latter uses the theory of reproducing kernel Hilbert spaces \cite{aronszajn1950, paulsen2016} and formulates system identification as function (operator) estimation in a possibly infinite-dimensional reproducing kernel Hilbert space. A key aspect of this paradigm is the \emph{representer} theorem, which allows for a computationally tractable solution to this estimation problem.
	
	A major challenge in system identification is the incorporation of prior knowledge about the physical properties of the system, such as stability \cite{pillonetto2018}, positivity \cite{grussler2017,khosravi2019} and (incremental) dissipativity \cite{rapisarda2011, vanwaarde2022}. The incorporation of such properties can enhance the reliability and robustness of the identified model by ensuring that the simulated behaviour aligns with the real-world behaviour of the system. This is particularly relevant when the model is used to simulate trajectories far away from the data set used for identification. For example, it has been shown \cite{umenberger2019, shakib2022} that incorporating stability in the identified model results in better long-term predictions and generalizations to new data sets.
	
	In this paper, we focus on the identification of passive systems, which can be seen as a special case of dissipative systems. Identification of dissipative systems has been studied in the case where the underlying model class consists of linear systems \cite{rapisarda2011}. On the other hand, identification of incrementally dissipative (nonlinear) systems has been studied in the case where the underlying model class is a reproducing kernel Hilbert space of operators \cite{vanwaarde2022}. The results in the latter are obtained by making use of the so-called \emph{scattering transform}, which allows to reformulate the original problem as the problem of identifying a system with a small incremental gain. While these results can be used to identify systems with non-incremental properties, i.e., systems that are simply dissipative, the reliance on the scattering transform presents some limitations. Most notably, the scattering transform is generally not computable in the non-incremental case.
	
	Motivated by this, we consider the problem of identifying a passive input-output model for a (nonlinear) system from data on its input and output trajectories. We consider an input-output model defined as an operator from the space of inputs to itself. In this case, passivity of the system corresponds to nonnegativity of the corresponding input-output operator, under the assumption that the latter is causal. Therefore, as a first step towards identifying a passive input-output model, we focus on the problem of identifying a nonnegative input-output operator. Since nonnegativity is not preserved under arbitrary linear combinations, we cannot directly apply the theory of reproducing Hilbert spaces and use the representer theorem. Instead, we take inspiration from the sum of squares idea for nonnegative function estimation recently developed in \cite{marteau-ferey2020}. 
	
	In particular, we consider a class of input-output operators defined via a feature map and a nonnegative bounded linear operator in such a way that nonnegativity of the input-output operator is guaranteed. Then, we formulate the identification problem as an estimation problem for the aforementioned linear operator.  As our main result, we obtain a representer theorem for this estimation problem in the case where the input space is finite-dimensional. Then, we show that a solution of the latter can be obtained by solving an associated semidefinite program. Finally, by means of an example, we demonstrate how these results can be used in practice, even when the input space is not finite-dimensional.
	
	The remainder of this paper is organized as follows. We introduce notation and preliminaries in Section~\ref{sec:notation}. We formulate the identification problem in Section~\ref{sec:problem} and we discuss the relevant theory of reproducing kernel Hilbert spaces in Section~\ref{sec:RKHS}. We solve the identification problem in Section~\ref{sec:nonnegative} and we show how the solution can be obtained by solving an associated semidefinite program in Section~\ref{sec:computation}. Finally, we demonstrate our results with an example in Section~\ref{sec:example} and we provide concluding remarks in Section~\ref{sec:conclusion}.

\section{Notation and prelimnaries}\label{sec:notation}

We denote the set of positive integers by $\bbN$. For $n\in\bbN$, we denote the set $\{1,\dots,n\}$ by $[n]$. We denote the Cartesian product of the sets $X$ and $Y$ by $X\times Y$. Given a set $X$ and $n\in\bbN$, we write $X^n$ to denote the Cartesian product
\begin{equation}
	X^n = \underbrace{X\times X \times \cdots \times X}_{n \text{ times}}.
\end{equation}
Given a Hilbert space $\calU$, we denote the corresponding inner product by $\inner{\cdot}{\cdot}_{\calU}$ and the induced norm by $\norm{\cdot}_{\calU}$. We say that an operator $G:\calU\to\calU$ is \emph{nonnegative} if $\inner{Gu}{u}_\calU \geq 0$ for all $u\in\calU$. Given another Hilbert space $\calW$, we denote the Banach space of bounded linear operators from $\calU$ to $\calW$ by $\B(\calU,\calW)$ and the operator norm by $\norm{\cdot}$. We use the shorthand notation $\B(\calU) = \B(\calU, \calU)$ and
\begin{equation}
	\B^+(\calU) = \set{G\in\B(\calU)}{G \text{ is nonnegative}}.
\end{equation}
We denote the adjoint of $B\in\B(\calU,\calW)$ by $B^* \in \B(\calW,\calU)$. We identify the vector {$w\in\calW$} with the operator $\alpha\mapsto\alpha w$, which belongs to $\B(\bbR,\calW)$. Consequently, $w^*\in\B(\calW,\bbR)$ is such that $w^*z = \inner{w}{z}_\calW$ for all $z\in\calW$. In particular, given $w_1,w_2\in\calW$, we have that $w_1w_2^*\in\B(\calW)$ is such that $w_1w_2^*z = \inner{w_2}{z}_\calW w_1$ for all $z\in\calW$.

Note that $\calU^n$ is a Hilbert space with inner-product
\begin{equation}
	\inner{(u_1,\dots,u_n)}{(v_1,\dots,v_n)}_{\calU^n} = \sum_{i=1}^{n}\inner{u_i}{v_i}_{\calU}.
\end{equation}
Given operators $M_{ij}\in\B(\calU)$, $i,j\in[n]$, we define
\begin{equation}\label{eq:block_operator}
	M = \bbm M_{11} & \cdots & M_{1n} \\ \vdots & \ddots & \vdots \\ M_{n1} & \cdots & M_{nn} \ebm
\end{equation}
as the operator in $\B(\calU^n)$ such that
\begin{equation}
	M(u_1,\dots,u_n) = \left(\sum_{i=1}^n M_{1i}u_i, \dots, \sum_{i=1}^n M_{ni}u_i\right).
\end{equation}
It is easily seen that any $M\in\B(\calU^n)$ can be written in the form \eqref{eq:block_operator} for suitably chosen  $M_{ij}\in\B(\calU)$, $i,j\in[n]$.

\subsection{Signal spaces}
\renewcommand{\d}{\textrm{d}}

Let $\bbR^+$ denote the set of nonnegative real numbers and let $\calR$ be a finite-dimensional Hilbert space, e.g., $\bbR^m$ with the Euclidian inner-product. The set of \emph{square integrable functions} $L_2(\bbR^+, \calR)$ consists of all measurable functions $u:\bbR^+\to\calR$ such that $\int_0^{\infty} \norm{u(t)}_\calR \d t < \infty$. If we identify functions that differ only on a set of Lebesgue measure zero, then $ L_2(\bbR^+, \calR)$ becomes a Hilbert space with inner-product
\begin{equation}\label{eq:square_integrable_inner}
	\inner{u}{v}_{L_2(\bbR^+, \calR)}= \int_{0}^{\infty} \inner{u(t)}{v(t)}_\calR \d t,
\end{equation}
The Hilbert space $L_2([0,T], \calR)$, where $T>0$, can be defined similarly by replacing the upper bound with $T$. 

Let $\bbZ^+$ denote the set of nonnegative integers. The set of \emph{square summable sequences} $\ell_2(\bbZ^+,\calR)$ consists of all functions $u:\bbZ^+\to\calR$ such that $\sum_{t=0}^{\infty} \norm{u(t)}_\calR< \infty.$ The latter is a Hilbert space with inner-product
\begin{equation}
	\inner{u}{v}_{\ell_2(\bbZ^+,\calR)} = \sum_{t=0}^{\infty} \inner{u(t)}{v(t)}_{\calR}.
\end{equation}
The Hilbert space $\ell_2(\{0,1\dots,T\}, \calR)$, where $T\in\bbN$, can be defined similarly.

We will collectively refer to the spaces above as \emph{signal spaces}. Thus, a signal space consists of functions  $u:\bbT\to\calR$, where $\calR$ is a Hilbert space and $\bbT$ is a time axis given by one of $\bbR^+$, $[0,T]$, $\bbZ^+$ or $\{0,1,\dots,T\}$. 

\subsection{Passivity and causality}
Next, we define the notions of causality and passivity, see \cite[Chapter~VI]{desoer2009} for details. Let $\calU$ be a signal space. For $\tau \in \bbT$, we define the truncation operator $P_\tau:\calU \to\calU$ as
\begin{equation}
	(P_\tau u)(t) = 
	\begin{dcases*}
		u(t) & $t \leq \tau$,\\
		0 & $t > \tau$.
	\end{dcases*}
\end{equation}
Note that $P_\tau$ is linear, $P_\tau^2 = P_\tau$ and $P_\tau^* = P_\tau$, that is, $P_\tau$ is an orthogonal projection for all $\tau\in\bbT$. We say that an operator $G:\calU\to\calU$ is:
\begin{enumerate}
	\item \emph{passive} if
	\begin{equation}\label{eq:passivity}
		\inner{P_\tau u}{P_\tau Gu}_\calU \geq 0
	\end{equation}
	for all $u\in\calU$ and $\tau\in\bbT$;
	\item \emph{causal} if 
	\begin{equation}
		P_\tau Gu = P_\tau G P_\tau u,
	\end{equation}
	for all $u\in\calU$ and $\tau\in\bbT$.
\end{enumerate}
In other words, $G$ is causal if the output at each time instant depends only on the past of the input. Note that $G$ is passive only if it is nonnegative (take $\tau\to\infty$ in \eqref{eq:passivity}). If $G$ is causal, then $G$ is passive if and only if $G$ is nonnegative since
\begin{equation*}
	\inner{P_\tau u}{P_\tau Gu}_\calU = \inner{P_\tau u}{P_\tau GP_\tau u}_\calU = \inner{P_\tau u}{ GP_\tau u}_\calU,
\end{equation*}
for all $u\in\calU$ and $\tau\in\bbT$, where the second equality follows from the fact that $P_\tau$ is an orthogonal projection.

\section{Problem formulation}\label{sec:problem}
Let $\calU$ be a signal space. Suppose that we are given a data set consisting of $n\in\bbN$ input-output signals
\begin{equation}\label{eq:data}
	(u_i,y_i)\in\calU\times\calU,\quad i\in[n].
\end{equation}
At a high-level, we are interested in finding a passive operator $G:\calU\to\calU$ that fits the data well, i.e., such that the \emph{data misfit}
\begin{equation}\label{eq:misfit}
	L(G) = \sum_{i=1}^{n} \norm{Gu_i - y_i}_\calU^2
\end{equation}
is small. Note that a passive operator is necessarily nonnegative. Therefore, as a first step towards solving this  high-level problem, in this paper, we focus on the following problem.
\begin{problem}\label{prob:nonnegativity}
	Given the data in \eqref{eq:data}, find a nonnegative operator $G:\calU\to\calU$ such that $L(G)$ is small.
\end{problem}

We will formalize what ``small" means by introducing an appropriate cost functional in Section~\ref{sec:nonnegative}. Then, we will solve Problem~\ref{prob:nonnegativity} in the case where $\calU$ is finite-dimensional. Note that the only finite dimensional signal space that we consider is $\ell_2(\{0,\dots,T\}, \calR)$, where $T\in\bbN$ and $\calR$ is a Hilbert space. Nevertheless, since nonnegativity is defined for operators on arbitrary Hilbert spaces, the following results are equally applicable to any finite-dimensional Hilbert space $\calU$. In particular, they are applicable to any finite dimensional subspace of a signal space. We will demonstrate this with an example in Section~\ref{sec:example}.

\section{Reproducing kernel Hilbert spaces}\label{sec:RKHS}

In this section, we will introduce part of the theory on reproducing kernel Hilbert spaces that is relevant to this paper. We refer to \cite{paulsen2016} for a detailed treatment of the case where the kernel is scalar-valued, and to \cite{micchelli2005} for the extension to the general case. We also refer to the technical report by \cite{kalnishkan2009}. In the context of this paper, we restrict the discussion to the case where the input and output spaces are identical.

Consider a Hilbert space $\calU$. A Hilbert space $\calH$ of operators from $\calU$ to $\calU$ is a \emph{reproducing kernel Hilbert space} if it admits a reproducing kernel, defined below.
\begin{definition}
	A map $k:\calU\times\calU\to\B(\calU)$ is a \emph{reproducing kernel} for $\calH$ if $k(\cdot, u)v:\calU \to\calU$ is a member of $\calH$ for all $u,v\in\calU$, and
		\begin{equation}\label{eq:reproducing_property}
			\inner{v}{Hu}_\calU = \inner{H}{k(\cdot, u)v}_\calH
		\end{equation}
		for all $u,v\in\calU$ and $H\in\calH$.
\end{definition} 

Property \eqref{eq:reproducing_property} is referred to as the \emph{reproducing property}. The class of reproducing kernels is completely characterized by two properties, namely, symmetry and positive semidefiniteness, defined below.
\begin{definition}
	The map $k:\calU\times\calU\to\B(\calU)$ is \emph{symmetric} if $k(u,v)^* = k(v,u)$ for all $u,v\in\calU,$.
\end{definition}
\begin{definition}
	The map $k:\calU\times\calU\to\B(\calU)$ is \emph{positive semidefinite} if 
	\begin{equation}\label{eq:pos_ker}
		\sum_{i=1}^n \inner{v_i}{k(u_i,u_j)v_j} \geq 0
	\end{equation}
	for all $n\in\bbN$, $u_1,\dots,u_n\in\calU$, $v_1,\dots,v_n\in\calU$.
\end{definition}

It is straightforward to show that a reproducing kernel needs to be symmetric and positive semidefinite due to the reproducing property. More notably, the converse also holds and is known as the Moore-Aronszajn theorem \cite{aronszajn1950}, see \cite{micchelli2005} for the extension to the general case.

\subsection{Gram operator}

The properties of symmetry and positive semidefiniteness of the map $k:\calU\times\calU \to\B(\calU)$ can be expressed using the associated \emph{Gram operator}. Given $n\in\bbN$ and $u_1,\dots,u_n\in\calU$, we define the Gram operator $K\in\B(\calU^n)$ as
\begin{equation}\label{eq:Gram}
	K = \bbm k(u_1,u_1) & \cdots & k(u_1,u_n) \\ \vdots & \ddots & \vdots \\ k(u_n,u_1) & \cdots & k(u_n,u_n) \ebm,
\end{equation}
Then, \eqref{eq:pos_ker} can be written as
\begin{equation}
	\inner{K(v_1,\dots,v_n)}{(v_1,\dots,v_n)}_{\calU_n} \geq 0.
\end{equation}
Therefore, the map $k$ is symmetric and positive semidefinite if and only if, for all $n\in\bbN$ and $u_1,\dots,u_n\in\calU$, the corresponding Gram operator $K$ in \eqref{eq:Gram} is self-adjoint and nonnegative.

\subsection{Feature map}
Symmetric positive definite maps  $k:\calU\times\calU\to\B(\calU)$ can be defined via an auxiliary Hilbert space $\calW$ and a map $\phi: \calU\to\B(\calU,\calW)$. In particular, it can be shown that the map $k(u,v) = \phi(u)^*\phi(v)$ is symmetric and positive semidefinite. The space $\calW$ is typically referred to as a \emph{feature space} and the map $\phi$ as a \emph{feature map}. There is a close relationship between symmetric positive definite maps $k$ and their associated feature maps $\phi$. This is captured in the following proposition.
\begin{proposition}
	The map $k:\calU\times\calU\to\B(\calU)$ is symmetric and positive semidefinite if and only if there exists a Hilbert space $\calW$ and a feature map $\phi: \calU\to\B(\calU,\calW)$ such that $k(u,v) = \phi(u)^*\phi(v)$.
\end{proposition}

The feature map $\phi$ is central in developing kernel-based algorithms and can allow for simpler analysis. In particular, it can be shown that the reproducing kernel Hilbert space $\calH$ with reproducing kernel $k(u,v) =  \phi(u)^*\phi(v)$ consists of operators of the form
\begin{equation}
	Hu = \phi(u)^* w,
\end{equation}
where $w\in\bar \calW$ and $\bar\calW\subset\calW$ is the closed linear span of the set $\set{\phi(v_1)v_2}{v_1,v_2\in\calU}$. We refer to \cite[Theorem~5 and Theorem~6]{kalnishkan2009} for a proof in the case where the kernel is scalar-valued. The proof in the general case follows similar reasoning.

\subsection{Representer theorem}

One of the most attractive features of the theory of reproducing kernel Hilbert spaces is that operator estimation problems have computationally tractable solutions. Here, we focus on the problem stated in Section~\ref{sec:problem} but without the nonnegativity requirement. In particular, we consider the problem of identifying an operator $G:\calU\to\calU$ such that the data misfit \eqref{eq:misfit} is small. Let $\calH$ be a reproducing kernel Hilbert space of operators from $\calU$ to $\calU$. We restrict our search to operators $G$ that belong to $\calH$. In general, under mild conditions, the data $(u_i,y_i)\in\calU\times\calU$, $i\in[n]$, can be interpolated exactly, i.e., there exists an operator $G\in\calH$ such that $Gu_i = y_i$ for all $i\in[n]$ and, thus, $L(G) = 0$. However, we typically expect the data to be noisy, in which case exact interpolation is undesirable. In such cases, a well-established identification approach is to instead consider the \emph{regularized least squares} problem
\begin{equation}\label{eq:min}
	\min_{G\in\calH}\quad L(G) +\gamma \norm{G}_\calH.
\end{equation}
Here, the first term promotes small $L(G)$, i.e., a good fit of the data, whereas the second term promotes small $\norm{G}_\calH$, which is interpreted as low complexity. The parameter $\gamma > 0$ regulates the balance between these two terms. 

The representer theorem, stated in the following proposition, shows that \eqref{eq:min} has a unique solution for any reproducing Hilbert space $\calH$, any parameter $\gamma >0$ and any data set $(u_i,y_i)\in\calU\times\calU$, $i\in[n]$.

\begin{proposition}\label{prop:representer}
	There exists a unique solution $\hat G$ of \eqref{eq:min}, which is given by
	\begin{equation}\label{eq:Ghat_k}
		\hat Gu = \sum_{i=1}^{n}k(u, u_i)c_i,
	\end{equation}
	with $c_i\in\calU$, $i\in[n]$, given by the unique solution of
	\begin{equation}\label{eq:representer_coeff}
		(K+\gamma I)(c_1,\dots,c_n) = (y_1,\dots,y_n),
	\end{equation}
	where $K\in\B(\calU^n)$ is the Gram operator associated with $k$ and $u_1,\dots,u_n$, and $I\in\B(\calU^n)$ is the identity operator.
\end{proposition}

Note that, if $\calU$ is finite-dimensional, then $\calU^n$ is finite-dimensional and, thus, solving \eqref{eq:representer_coeff} amounts to solving a finite-dimensional system of linear equations.

As already mentioned in the previous subsection, there exists feature map $\phi:\calU\to\B(\calU,\calW)$ such that the reproducing kernel for $\calH$ is given by $k(u,v) = \phi(u)^*\phi(v)$. Furthermore, every $G\in\calH$ is of the form
\begin{equation}
	Gu = \phi(u)^* w,
\end{equation}
where $w\in \bar \calW\subset\calW$. With a minor abuse of notation, we can treat the data misfit \eqref{eq:misfit} as a function of $w$, such that
\begin{equation}
	L(w) = \sum_{i=1}^n \norm{\phi(u_i)^* w - y_i}_\calU^2.
\end{equation}
Since $\norm{G}_\calH = \norm{w}_\calW$, we can rewrite \eqref{eq:min} as
\begin{equation}\label{eq:min_w}
	\min_{w\in\bar\calW}\quad L(w) +\gamma \norm{w}_\calW.
\end{equation}
Then, Proposition~\ref{prop:representer} tells us that \eqref{eq:min_w} has a unique solution $\hat w \in\bar \calW$ given by
\begin{equation}\label{eq:w_representer}
	\hat w = \sum_{i=1}^n \phi(u_i)c_i,
\end{equation}
with $c_i\in\calU$, $i\in[n]$, given by the unique solution of \eqref{eq:representer_coeff}.

\section{Representer theorem for nonnegative operators}\label{sec:nonnegative}
In this section, we will provide a solution to Problem~\ref{prob:nonnegativity}. To do this, we take inspiration from the theory of reproducing kernel Hilbert spaces and the sum of squares idea developed in \cite{marteau-ferey2020}. We will first motivate our approach by considering a simple example.

\begin{example}\label{ex:approach}
	Consider a continuous function $f:\bbR\to\bbR$ that is differentiable at zero. We view $f$ as an operator from $\bbR$ to $\bbR$. As such, $f$ is nonnegative if and only if $f(x)x \geq 0$ for all $x\in\bbR$. In particular, this means that $f(0) = 0$, hence we can write $f(x) = g(x)x$, where the function $g:\bbR\to\bbR$ given by
	\begin{equation}
		g(x) = \begin{dcases*}
			\frac{f(x)}{x} & if $x\neq 0$,\\
			f'(0)& if $x=0$.
		\end{dcases*}
	\end{equation}
	is continuous. Since $f(x)$ is nonnegative as an operator, it follows that $g(x) \geq 0$ for all $x\in\bbR$, i.e., $g$ is nonnegative as a function. Conversely, if $f(x) = g(x)x$ for some nonnegative function $g$, then $f$ is nonnegative as an operator. Now, suppose that we have a data set consisting of $n\in\bbN$ pairs $(x_i,f_i)\in\bbR\times\bbR$, $i\in[n]$, and consider the problem of finding a nonnegative operator $f:\bbR\to\bbR$ that fits the data well, i.e., $f(x_i) \approx f_i$ for all $i\in[n]$. If we further require that $f$ is continuous and differentiable at $0$, then, without loss of generality, we can assume that $f(x) = g(x)x$ for some continuous nonnegative function $g:\bbR\to\bbR$. Consequently, the problem reduces to finding a nonnegative $g$ such that $g(x_i)\approx f_i/x_i$ for all $i\in[n]$. This can be solved using the sum-of-squares approach in \cite{marteau-ferey2020}. In particular, given a Hilbert space of features $\calW$ and a feature map $\phi:\bbR\to\calW$, we can model a nonnegative function as $g(x) = \phi(x)^* Q \phi(x)$, where the operator $Q\in\B(\calW)$ is nonnegative. Note that this yields
	\begin{equation}\label{eq:f_form}
		f(x) = \phi(x)^* Q \phi(x)x
	\end{equation}
	As shown in \cite{marteau-ferey2020}, the problem can be reformulated as a regularized least squares problem with an appropriately chosen cost functional. Furthermore, the authors of \cite{marteau-ferey2020} obtain a representer theorem for this problem that provides a computationally tractable solution.
\end{example}

With Example~\ref{ex:approach} in mind, we take the following approach to solve Problem~\ref{prob:nonnegativity}. Let $\calW$ be a Hilbert space of features and $\phi: \calU\to\B(\calU,\calW)$ a feature map. In view of \eqref{eq:f_form}, we consider input-output operators $G:\calU\to\calU$ of the form
\begin{equation}\label{eq:G}
	Gu = \phi(u)^* Q \phi(u) u,
\end{equation}
where $Q\in\B(\calW)$.  Since
\begin{equation*}
	\inner{Gu}{u}_\calU = \inner{ \phi(u)^* Q \phi(u) u}{u}_\calU = \inner{Q\phi(u)u}{\phi(u)u}_\calW,
\end{equation*}
it follows that that $G$ is nonnegative if $Q$ is nonnegative. With a minor abuse of notation, we can treat the data misfit \eqref{eq:misfit} as a function of $Q$, such that
\begin{equation}
	L(Q) = \sum_{i=1}^n \norm{\phi(u_i)^* Q \phi(u_i) u_i- y_i}_\calU^2.
\end{equation}
We can now reformulate Problem~\ref{prob:nonnegativity} as the regularized least squares problem
\begin{equation}\label{eq:min_re}
	\min_{Q\in \B^+(\calW)} L(Q) +\gamma \norm{Q},
\end{equation}
where $\norm{\cdot}$ is the operator norm and $\gamma>0$. We can obtain a representer theorem for this problem under the condition that $\calU$ is finite-dimensional. Therefore, for the remainder of this paper, we will assume that $\calU$ is finite-dimensional.

Consider the subspace $\hat\calW\subset\calW$ given by
\begin{equation}
	\hat \calW = \sum_{i=1}^n \im \phi(u_i)
\end{equation}
Since $\calU$ is finite-dimensional, it follows that $\hat\calW$ is also finite-dimensional. Let $\Pi\in\B(\calW)$ be the orthogonal projection onto $\hat\calW$. This means that $\im \Pi = \hat\calW$, $\Pi^2 = \Pi$ and $\Pi^* = \Pi$.  Consider the following lemma.
\begin{lemma}\label{lem:projection_nonnegativity}
	The following statements hold:
	\begin{enumerate}
		\item $L(\Pi Q \Pi) = L(Q)$ for all $Q\in\B(\calW)$;
		\item $\norm{\Pi Q \Pi} \leq \norm{Q}$ for all $Q\in\B(\calW)$;
		\item $\Pi Q \Pi\in\B^+(\calW)$ for all $Q\in\B^+(\calW)$.
	\end{enumerate}
\end{lemma}\
\proof
	Statement 1) holds because $\Pi \phi(u_i) = \phi(u_i)$ and $\phi(u_i)^* \Pi = \phi(u_i)^*$ for all $i\in[n]$. We also have that
	\begin{equation}
		\inner{\Pi Q \Pi w}{w}_\calW = \inner{Q\Pi w}{\Pi w}_\calW
	\end{equation}
	for all $Q\in\B(\calW)$ and $w\in\calW$, hence statement 3) holds. To show that statement 2) holds as well, let $Q\in\B(\calW)$ and note that $\norm{\Pi w}_\calW\leq \norm{w}_\calW$ for all $w\in\calW$. Therefore, we have that $\norm{\Pi Q \Pi w}_\calW \leq \norm{Q \Pi w}_\calW$ and $\norm{w}_\calW \leq 1$ implies $\norm{\Pi w}_\calW \leq 1$ for all $w\in\calW$. This implies that
	\begin{align*}
		\norm{\Pi Q \Pi} &= \sup_{\norm{w}\leq 1} \norm{\Pi Q \Pi w}_\calW\\
		&\leq \sup_{\norm{w}\leq 1} \norm{Q \Pi w}_\calW\\
		&\leq \sup_{\norm{\Pi w}\leq 1} \norm{Q \Pi w}_\calW\\
		&\leq \norm{Q}
	\end{align*}
	and, thus, statement 3) holds.
\endproof

The first two statements of Lemma~\ref{lem:projection_nonnegativity} tell us that the cost in \eqref{eq:min_re} does not increase when $Q$ is replaced by $\Pi Q \Pi$. Together with the third statement, this implies that
\begin{equation}\label{eq:min_hatQ}
	\min_{ Q\in \B^+(\calW)} L( Q) +\gamma \norm{ Q} = \min_{\hat Q\in \hatQ} L(\hat Q) +\gamma \norm{\hat Q}, 
\end{equation}
where $\hatQ = \set{\Pi Q \Pi}{Q\in\B^+(\calW)} \subset\B^+(\calW)$. In the following, we will show that $\hatQ$ has a finite parametrization. To this end, let the operator $\Phi\in\B(\calU^n,\calW)$ be given by
\begin{equation}
	\Phi = \bbm \phi(u_1) & \cdots & \phi(u_n) \ebm
\end{equation}
and consider the following lemma.

\begin{lemma}\label{lem:projection_form}
	There exist a self-adjoint operator $N\in\B(\calU^n)$ such that $\Pi = \Phi N \Phi^*$.
\end{lemma}
\proof
Suppose that the dimension of $\hat\calW$ is $K\in\bbN$ and let $\{w_1,\dots,w_K\}$ be an orthonormal basis for $\hat\calW$. It is well-known \cite[Proposition~4.7]{conway2007} that $\Pi = \sum_{k=1}^{K} w_iw_i^*$. Since $w_k\in\hat \calW$ for all $k\in[K]$, there exist $v_{k,i}\in\calU$, $i\in\oneto{n}$, such that $w_k = \sum_{i=1}^{n} \phi(u_i)v_{k,i}$. Consequently, we obtain
\begin{equation}
	\Pi = \sum_{k=1}^K\sum_{i,j=1}^{n}  \phi(u_i)v_{k,i} v_{k,j}^*\phi(u_j)^*
\end{equation}
hence the operators $N_{ij}= \sum_{k=1}^{K} v_{k,i} v_{k,j}^* $, $i,j\in[n]$, are such that $\Pi = \sum_{i,j=1}^n \phi(u_i) N_{ij} \phi(u_j)^*$. In other words, we have that $\Pi = \Phi N \Phi^*$ with $N\in\B(\calU^n)$ given by
\begin{equation}
	N = \bbm N_{11} & \cdots & N_{1n} \\ \vdots & \ddots & \vdots \\ N_{n1} & \cdots & N_{nn} \ebm
\end{equation}
Since $N_{ij}^* = N_{ji}$ for all $i,j\in[n]$, it follows that $N$ is self-adjoint, which concludes the proof.
\endproof
The following lemma is a direct consequence of Lemma~\ref{lem:projection_form}.
\begin{lemma}\label{lem:Q_form}
	We have that $\hatQ = \set{\Phi M \Phi^*}{M\in\B^+(\calU^n)}$.
\end{lemma}
\proof
Let $\hat Q = \Phi M \Phi^*$ for some $M\in\B^+(\calU^n)$. Since $\Pi \phi(u_i) = \phi(u_i)$ for all $i\in[n]$, it follows that $\hat Q = \Pi \hat Q \Pi$. Furthermore, since $M$ is nonnegative, it follows that 
\begin{equation*}
	\inner{\hat Qw}{w}_\calW = \inner{\Phi M\Phi^* w}{ w}_\calW = \inner{M\Phi^* w}{\Phi^* w}_{\calU^n} \geq 0
\end{equation*}
for all $w\in\calW$, hence $\hat Q \in \B^+(\calW)$ and, thus, $\hat Q\in\hatQ$. For the converse, let $\hat Q \in \hatQ$, i.e., $\hat Q = \Pi Q \Pi$ with $Q\in\B^+(\calW)$. Due to Lemma~\ref{lem:projection_form}, there exists a self-adjoint $N\in\B(\calU^n)$ such that $\Pi = \Phi N \Phi^*$, hence $\hat Q = \Phi M \Phi^*$ with $M\in\B(\calU^n)$ given by $M = N \Phi^* Q \Phi N$. As $Q$ is nonnegative and $N$ is self-adjoint, it follows that $M\in\B^+(\calU^n)$,
\endproof

Note that $\B(\calU^n)$ is finite-dimensional since $\calU$ is finite-dimensional. Therefore, Lemma~\ref{lem:Q_form} shows that $\hatQ$ has a finite parametrization. Furthermore, due to Lemma~\ref{lem:projection_nonnegativity}, Lemma~\ref{lem:Q_form} and \eqref{eq:min_hatQ}, a solution of \eqref{eq:min} can be obtained by solving
\begin{equation}\label{eq:min_M}
	\min_{ M \in \B^+(\calU^n)} L(\Phi M \Phi^*) + \gamma \norm{\Phi M \Phi^*}.
\end{equation}
With this in mind, we obtain the following representer theorem for the problem \eqref{eq:min}, which is the main result of this paper.
\begin{theorem}\label{thm:representer}
	Suppose that $\calU$ is finite-dimensional. A solution of \eqref{eq:min_re} exists if and only if a solution of \eqref{eq:min_M} exists. Moreover, if it exists, a solution of \eqref{eq:min_re} is given by
	\begin{equation}\label{eq:Q_representer_compact}
		\hat Q = \Phi M \Phi^*,
	\end{equation}
	where $M\in\B^+(\calU^n)$ is a solution of \eqref{eq:min_M}.
\end{theorem}
\proof
Suppose that $Q\in\B^+(\calW)$ is a solution of \eqref{eq:min_re}. Due to Lemma~\ref{lem:projection_nonnegativity}, it follows that $\Pi Q \Pi$ is nonnegative and
\begin{equation}\label{eq:PQP_Q}
	L(\Pi  Q \Pi) + \norm{\Pi  Q \Pi} \leq L( Q) + \norm{ Q}.
\end{equation}
Since $Q$ is a solution of  \eqref{eq:min_re}, the latter holds with equality and $\hat Q = \Pi Q \Pi \in \hatQ$ is also a solution of \eqref{eq:min_re}. Due to Lemma~\ref{lem:Q_form}, there exists $M\in\B^+(\calU^n)$ such that $\hat Q = \Phi M \Phi^*$. As $\hat Q$ is a solution of \eqref{eq:min_re}, we obtain
\begin{equation}
	L(\hat Q) + \norm{\hat Q} \leq L(Q') + \norm{Q'}
\end{equation}
for all $Q' \in\B^+(\calW)$ and, in particular, for all $Q' = \Phi M' \Phi^*$, where $M'\in\B^+(\calU^n)$. This implies that 
\begin{equation}\label{eq:M_M'}
	L(\Phi M \Phi^*) + \norm{\Phi M \Phi^*} \leq L(\Phi M' \Phi^*) + \norm{\Phi M' \Phi^*}
\end{equation}
for all $M'\in\B^+(\calU^n)$ and, thus, $M$ is a solution of \eqref{eq:min_M}.

To show the converse, suppose that  $M$ is a solution of \eqref{eq:min_M}. We will show that $\Phi M \Phi^* \in \B^+(\calW)$ is a solution of \eqref{eq:min_re}. Let $Q\in\B^+(\calW)$ and note that \eqref{eq:PQP_Q} holds due to Lemma~\ref{lem:projection_nonnegativity}. Furthermore, due to Lemma~\ref{lem:Q_form}, there exists $M'\in\B^+(\calU^n)$ such that $\Pi Q \Pi = \Phi M' \Phi^*$. Now, \eqref{eq:M_M'} holds because $M$ is a solution of \eqref{eq:min_M}. Together with \eqref{eq:PQP_Q}, this implies that
\begin{equation}
	L(\Phi M \Phi^*) + \norm{\Phi M \Phi^*} \leq L(Q) + \norm{ Q}
\end{equation}
Since the latter holds for arbitrary $Q\in\B^+(\calW)$, it follows that $\hat Q = \Phi M \Phi^*$ is a solution of \eqref{eq:min_re}.
\endproof

Recall that any $M\in\B(\calU^n)$ is of the form \eqref{eq:block_operator}, where $M_{ij}\in\B(\calU)$, $i,j\in[n]$. Therefore, we can rewrite \eqref{eq:Q_representer_compact} as
\begin{equation}\label{eq:Q_representer}
	\hat Q = \sum_{i,j=1}^{n} \phi(u_i) M_{ij} \phi(u_j)^*
\end{equation}
Note the parallel between \eqref{eq:min_w} and \eqref{eq:min_re}, as well as \eqref{eq:w_representer} and \eqref{eq:Q_representer}. In particular, while \eqref{eq:w_representer} is a linear combination of features, \eqref{eq:Q_representer} can be interpreted as a sum of squares of features. We point out that, unlike the representer theorem stated in Proposition~\ref{prop:representer}, Theorem~\ref{thm:representer} does not tell us that a solution of \eqref{eq:min_re} exists or that the solution is unique if it does exist. Nevertheless, Theorem~\ref{thm:representer} paves the way to a computationally tractable solution of Problem~\ref{prob:nonnegativity} in the case where $\calU$ is finite-dimensional.  Indeed, in the next section we will show that a solution of Problem~\ref{prob:nonnegativity} can be obtained by solving an associated semidefinite program. 

Before we do that, we note that the operator $G:\calU\to\calU$ corresponding to \eqref{eq:Q_representer_compact} is given by
\begin{equation}\label{eq:G_M}
	Gu = \kappa(u)^* M \kappa(u) u,
\end{equation}
where $\kappa(u)\in\B(\calU,\calU^n)$ is given by
\begin{equation}\label{eq:Ku}
	\kappa(u) = \Phi^*\phi(u) =  \bbm k(u_1, u) \\ \vdots \\ k(u_n, u) \ebm,
\end{equation}
and $k:\calU\times\calU\to\B(\calU)$ is the reproducing kernel corresponding to the the feature map $\phi$, i.e., $k(v_1,v_2) = \phi(v_1)^*\phi(v_2)$. This means that evaluation of $G$ does not require evaluation of the feature map $\phi$, but only of the kernel $k$ corresponding to $\phi$. In fact, it can be shown that solving \eqref{eq:min_M} also does not require evaluation of $\phi$.  To this end, note that
\begin{equation}\label{eq:cost_kernel}
	L(\Phi M \Phi^*) = \sum_{i=1}^{n} \norm{K_i^* M K_iu_i - y_i}^2_{\calU},
\end{equation}
where $K_i= \kappa(u_i)\in\B(\calU^n,\calU)$ for all $i\in[n]$. The Gram operator associated with $k$ and $u_1,\dots,u_n,$ is given by
\begin{equation}
	K = \bbm K_1 & \cdots & K_n \ebm = \Phi^*\Phi.
\end{equation}
Since $K\in\B(\calU^n)$ is self-adjoint and nonnegative, there exists a self-adjoint and nonnegative $K^\half\in\B(\calU^n) $ such that $K^\half K^\half = K$, and we have the following lemma.

\begin{lemma}
	We have that
	\begin{equation}\label{eq:PhiMPhi_norm}
		\norm{\Phi M \Phi^*} =  \Vert K^\half M K^\half \Vert
	\end{equation}
	for all $M\in \B(\calU^n)$.
\end{lemma}
\proof
As $\calU$ is finite-dimensional and $\Phi\in\B(\calU^n, \calW)$, it follows that $\im \Phi$ is finite-dimensional. Consequently, for all $w\in\calW$, we can write $w = \Phi x + z$, where $x\in\calU^n$ and $z \in (\im \Phi)^\perp = \ker \Phi^*$. Using this, we obtain
\begin{equation}
	\norm{\Phi M \Phi^* w}^2_\calW = \norm{\Phi M K x}^2_{\calW},
\end{equation}
where we used the fact that $K = \Phi^*\Phi$. This implies that
\begin{equation}
	\norm{\Phi M \Phi^*}^2 = \sup_{\norm{\Phi x}_{\calW} \leq 1} \norm{\Phi M K x}^2_{\calW}.
\end{equation}
Since $K^\half K^\half = \Phi^*\Phi$ and $K^\half$ is self-adjoint, we have that
\begin{equation}
	\norm{\Phi x}_{\calW}^2  = \inner{\Phi x}{\Phi x}_{\calW} = \Vert{K^\half x}\Vert_{\calU^n}^2,
\end{equation}
We also have that $\norm{\Phi M K x}^2_{\calW}= \Vert{K^\half M K x}\Vert_{\calU^n}^2$ and, thus,
\begin{equation}\label{eq:PhiMPhi_x}
	\norm{\Phi M \Phi^*}^2 = \sup_{\Vert{K^\half x}\Vert_{\calU^n} \leq 1} \Vert{K^\half M K x}\Vert^2_{\calU^n}.
\end{equation}
Since $\im K^\half$ is finite-dimensional, $z\in\calU^n$ can be written as $z = K^\half x + z'$ with $x\in\calU^n$ and $z' \in \im (K^\half)^\perp = \ker K^\half$. Then, \eqref{eq:PhiMPhi_x} implies that
\begin{equation*}
	\norm{\Phi M \Phi^*}^2 = \sup_{\Vert{z}\Vert_{\calU^n} \leq 1} \Vert{K^\half M K^\half z}\Vert^2_{\calU^n} = \Vert{K^\half M K^\half}\Vert^2,
\end{equation*}
and, thus, \eqref{eq:PhiMPhi_norm} holds.
\endproof

As an immediate consequence of \eqref{eq:cost_kernel} and \eqref{eq:PhiMPhi_norm}, we obtain the following corollary of Theorem~\ref{thm:representer}.
\begin{corollary}\label{cor:representer}
	If it exists, a solution of \eqref{eq:min_re} is given by \eqref{eq:Q_representer_compact}, where $M\in\B(\calU^n)$ is a solution of
	\begin{equation}\label{eq:min_M_re}
		\min_{M\in\B^+(\calU^n)} \sum_{i=1}^{n} \norm{K_i^* M K_iu_i - y_i}^2_\calU + \norm{K^\half M K^\half}.
	\end{equation}
	Moreover, the operator $G:\calU\to\calU$ corresponding to this solution is given by \eqref{eq:G_M}.
\end{corollary}

Corollary~\ref{cor:representer} shows that a solution of Problem~\ref{prob:nonnegativity} can be obtained by solving \eqref{eq:min_M_re}. Note that all operators in \eqref{eq:min_M_re} belong to $\B(\calU^n)$, which is finite-dimensional, hence Corollary~\ref{cor:representer} provides a computationally tractable solution of Problem~\ref{prob:nonnegativity}. In the next section, we will show how this solution can be obtained by solving an associated semidefinite program.

\section{Computional aspects}\label{sec:computation}

In this section, we will show how to obtain a solution of Problem~\ref{prob:nonnegativity} by solving an associated semidefinite program. To this end, since $\calU$ is finite-dimensional, it is isometrically isomorphic to $\bbR^m$ for some $m\in\bbN$. Therefore without loss of generality, we can assume that $\calU = \bbR^m$ with the Euclidian inner-product $\inner{\cdot}{\cdot}_2$ and induced norm $\norm{\cdot}_2$. Now, consider \eqref{eq:min_M_re}. Note that
$\B^+(\calU^n)$ consists of matrices $M\in\bbR^{nm\times nm}$ such that $M + M\t$ is symmetric positive semidefinite, which we denote by $M + M\t \geq 0$. Furthermore, we have that
\begin{equation}
	K_i\in\bbR^{mn\times m},\quad K^\half \in\bbR^{nm\times nm},\quad K^\half\in\bbR^{nm\times nm}
\end{equation}
where $K\geq 0$ and $K^\half \geq 0$. It is well known, see, e.g., \cite{vandenberghe1996}, that $\norm{K_i\t M K_iu_i - y_i}^2_2 \leq p_i$, $i\in[n]$, if and only if
\begin{equation}\label{eq:SDP_pi}
	\bbm I & K_i\t M K_iu_i - y_i \\ u_i\t K_i\t M\t K_i - y_i\t & p_i \ebm \geq 0.
\end{equation}
Similarly, $\norm{K^\half M K^\half} \leq p_0$ if and only if
\begin{equation}\label{eq:SDP_p0}
	\bbm p_0 I & K^\half M K^\half \\ K^\half M\t K^\half & p_0 I \ebm \geq 0.
\end{equation}
It can be shown that the latter is equivalent to
\begin{equation}\label{eq:SDP_p0_re}
	\bbm p_0K & KM K \\ KM\t K & p_0K \ebm \geq 0,
\end{equation}
whose evaluation does not require the computation of the matrix $K^\half$. Indeed, \eqref{eq:SDP_p0} holds if and only if
\begin{equation}
	p_0(x_1\t x_1 + x_2\t x_2) + 2x_1\t K^\half M K^\half x_2 \geq 0
\end{equation} 
for all $x_1,x_2\in\bbR^{mn}$. Since every $x\in\bbR^{mn}$ can be written as $x =  K^\half \bar x + z$, where $\bar x\in\bbR^{mn}$ and $z\in\ker K^{\half}$, it follows that \eqref{eq:SDP_p0} holds if and only if
\begin{equation*}
	p_0(\bar x_1\t K\bar x_1 + \bar x_2\t K\bar x_2 + z_1\t z_1 + z_2\t z_2) + 2\bar x_1\t K M K \bar x_2 \geq 0
\end{equation*}
for all $\bar x_1,\bar x_2\in\bbR^{mn}$ and $z_1,z_2\in\ker K^\half$. Since $z_1\t z_1 \geq 0$ and $z_2\t z_2 \geq 0$ for all $z_1,z_2\in\ker K^\half$, it follows that \eqref{eq:SDP_p0} holds if and only if 
\begin{equation}
	p_0(\bar x_1\t K\bar x_1 + \bar x_2\t K\bar x_2 ) + 2\bar x_1\t K M K \bar x_2 \geq 0,
\end{equation}
for all $\bar x_1,\bar x_2\in\bbR^{mn}$, that is, \eqref{eq:SDP_p0_re} holds.
With this in mind, the following corollary of Theorem~\ref{thm:representer} is an immediate consequence of Corollary~\ref{cor:representer}.

\begin{corollary}\label{cor:semi}
	Suppose that $\calU = \bbR^m$. If it exists, a solution of \eqref{eq:min_M_re} can be obtained by solving the semidefinite program
	\begin{alignat*}{2}
		\min& \quad &&\sum_{i=1}^n p_i + \gamma p_0 \label{eq:SDP}\\
		\text{s.t.}& \quad&&M + M\t \geq 0 \text{ and } \eqref{eq:SDP_pi}, \eqref{eq:SDP_p0_re} \text{ hold for all } i\in[n]
	\end{alignat*}
	Moreover, the operator $G:\bbR^m\to\bbR^m$ corresponding to this solution is given by \eqref{eq:G_M}.
\end{corollary}

\section{Illustrative example}\label{sec:example}

In this section, we will demonstrate how our results can be used in practice. Consider the nonlinear system
\begin{equation}\label{eq:ex_system}
	\begin{aligned}
		2\ddot{q} + \ddot{\theta} \cos(\theta) &=  \dot \theta^2\sin(\theta) - q\\
		2\ddot{\theta} + \ddot{q}\cos(\theta) &= u\\
		y &= \dot \theta
	\end{aligned}
\end{equation}
The latter describes the motion of a rotational/translational proof mass actuator \cite{bupp1995}. For simplicity, we have neglected disturbances and taken all parameters to be equal to 1. 

Let $T>0$ and $\calU = L_2([0,T], \bbR)$. We can view \eqref{eq:ex_system} as an operator $\Gtrue:\calU\to \calU$ that maps the input signal $u\in \calU$ to the corresponding output signal $\Gtrue u\in \calU$ generated by \eqref{eq:ex_system}, where the system is initially at rest, i.e.,
\begin{equation}
	q(0) = 0,\quad \dot{q}(0) = 0,\quad \theta(0) = 0,\quad \dot{\theta}(0) = 0.
\end{equation}
It can be shown that $\Gtrue$ is passive, hence, in particular, $\Gtrue$ is nonnegative. Now, suppose that we are given a data set consisting of $n\in\bbN$ input-output signals
\begin{equation}\label{eq:ex_data_traj}
	(u_i,y_i)\in \calU\times \calU,\quad i\in[n],
\end{equation}
generated by the system \eqref{eq:ex_system}. We are interested in identifying a nonnegative operator $G: \calU\to \calU$ that fits the data well, i.e., $Gu_i \approx y_i$ for all $i\in[n]$. Since $\calU$ is infinite-dimensional, we cannot use our results directly. We circumvent this issue by restricting to a finite-dimensional subspace of $\calU$.

To this end, we can approximate signals in $\calU$ to arbitrary precision  by a finite-degree polynomial. In particular, we can approximate signals using the first $m$ Legendre polynomials \cite{canuto2006}, where $m$ is as large as necessary to achieve a desired accuracy of the approximation. We refer the reader to \cite{rapisarda2023} for an example of the use of Legendre polynomials in control. 

Recall \cite[Section~2.3.2]{canuto2006} that the Legendre polynomials on the interval $[-1, 1]$ are given by $L_1(t) = 1$, $L_2(t) = t$, and the recurrence relation
\begin{equation}
	\frac{\d}{\d t}L_{i+2} = \frac{\d}{\d t}L_{i} + (2i+2)L_{i+1},\quad i\in\bbN.
\end{equation}
The Legendre polynomials on the interval $[0,T]$ can be obtained via the transformation $t\mapsto \frac{T}{2}(t+1)$. We assume that each Legendre polynomial is scaled so that $\norm{L_i}_\calU = 1$ for all $i\in\bbN$. As such, the Legendre polynomials form an orthonormal basis for $\calU$, i.e., every $u\in \calU$ can be written as $u = \sum_{i=1}^\infty c_i L_{i}$ for some coefficients $c_i\in\bbR$, $i\in\bbN$, which we call the \emph{Legendre coefficients} of $u$. Therefore, $\Gtrue$ is associated with an operator $\tilde{G}_{\textup{true}}$ that maps the Legendre coefficients of $u$  to the Legendre coefficients of $\Gtrue u$. Since the Legendre polynomials are orthonormal, $\Gtrue$ is nonnegative if and only if $\tilde{G}_{\textup{true}}$ is nonnegative.

With this in mind, let $m\in\bbN$ be such that the first $m$ Legendre polynomials can be used to accurately approximate the data set \eqref{eq:ex_data_traj}. Given a desired precision, an appropriate $m$ can be determined by using existing error bounds, see \cite[Section~5.4]{canuto2006} for details. The approximation yields a data set consisting of $n\in\bbN$ coefficient vectors
\begin{equation}\label{eq:ex_data_coefficients}
	(\tilde u_i, \tilde y_i)\in\bbR^m\times\bbR^m,\quad i\in[n],
\end{equation}
such that $u_i \approx \sum_{j=1}^m\tilde u_{ij}L_i$ and $y_i \approx \sum_{j=1}^m \tilde y_{ij} L_{i}$, where $\tilde u_{ij}$ and $\tilde y_{ij}$ are the $j$'th entries of $\tilde u_i$ and $\tilde y_i$, respectively. Now, we are interested in identifying a nonnegative operator $\tilde G:\bbR^m\to\bbR^m$ that fits the data well, i.e., $\tilde G\tilde u_i \approx \tilde y_i$ for all $i\in[n]$. This would correspond to a nonnegative operator $G:\tilde\calU \to\tilde\calU$, where $\tilde\calU\subset\calU$ is given by the linear span of the first $m$ Legendre polynomials. This operator $G$ is expected to be a good approximation of $\Gtrue$ when the input signal is restricted (or close) to $\tilde\calU$.

\begin{figure}
	\color{black}
	\centering
	\input{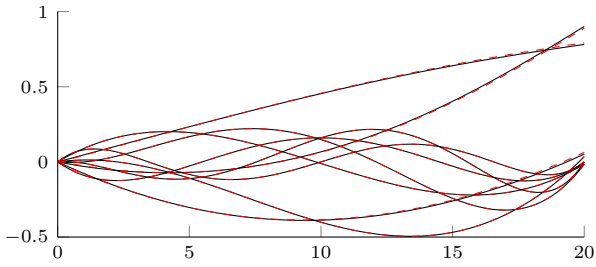}
	\vspace{-2mm}
	\caption{True (black, solid) and estimated data outputs with the Gaussian kernel (red, dashed).}
	\label{fig:ex_Gauss_data_fit}
	\vspace{-3mm}
\end{figure}

To illustrate the effectiveness of this approach, we consider the following specific example. Let $T = 20$, $n=9$ and consider the data set \eqref{eq:ex_data_traj}, where
\begin{equation}
	u_i = \begin{dcases}
		L_{i} & i\in\{1,\dots, 5\},\\
		L_{i-5} + L_{i-4} & i\in\{6,\dots, 9\}.
	\end{dcases}
\end{equation}
This means that the first 5 Legendre polynomials can be used to exactly represent the input signals $u_i$, $i\in[9]$. Approximating the corresponding output signals $y_i$, $i\in[9]$, using the first 10 Legendre polynomials yields relative errors less than $10^{-4}$. This is small enough for the purposes of this example, hence we take $m=10$. Next, we solve the semidefinite program in Corollary~\ref{cor:semi}, where $\gamma = 10^{-3}$ and the data set is given by the corresponding Legendre coefficient vectors \eqref{eq:ex_data_coefficients}. Furthermore, we use a simple Gaussian kernel 
\begin{equation}
	k(\tilde v_1,\tilde v_2) = e^{\frac{-\norm{\tilde v_1 - \tilde v_2}^2_2}{100^2}} I,
\end{equation}
where $\norm{\cdot}_2$ is the Euclidian norm and $I\in\bbR^{10\times 10}$ is the identity matrix. As can be seen on Figure~\ref{fig:ex_Gauss_data_fit}, the identified operator fits the data quite well; the true and estimated outputs are almost indistinguishable and the data misfit is approximately $1.6547\cdot10^{-4}$. To test the identified operator, we consider 1000 input signals given by linear combinations of the first 5 Legendre polynomials with coefficients from a uniform distribution on the interval $[0,1]$. Then, we compute the relative error between the true output signal generated by \eqref{eq:ex_system} and the estimated output signal generated by the identified operator. We find that the average relative error is approximately $3.7949\cdot 10^{-2}$. 
\begin{figure}
	\color{black}
	\centering
	\input{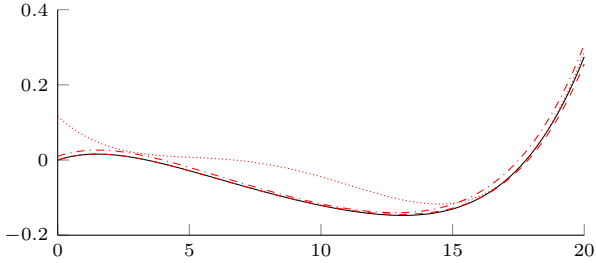}
	\vspace{-2mm}
	\caption{True (black, solid) and estimated outputs with the Gaussian (red, dashed), Laplacian (red, dot-dashed) and bilinear (red, dotted) kernels.}
	\label{fig:comparison}
	\vspace{-5mm}
\end{figure}

We perform the same tests with the identified operators when using the Laplacian and bilinear kernels, given by
\begin{equation}
	k(\tilde v_1, \tilde v_2) = e^{\frac{-\norm{\tilde v_1 - \tilde v_2}_2}{100}} I \quad\text{and}\quad k(\tilde v_1,\tilde v_2) = \tilde v_1\t \tilde v_2I,
\end{equation}
respectively. We find that the data misfit and average relative error when using the Laplacian kernel are approximately $2.2123\cdot10^{-7}$ and $8.3060\cdot10^{-2}$, respectively. On the other hand, the data misfit and average relative error when using the bilinear kernel are approximately $3.5530\cdot10^{-7}$ and 1.3889, respectively. An example comparison between the true output signal and the estimated output signals of each of these identified operators is seen on Figure~\ref{fig:comparison}. We see that the Gaussian kernel performs the best, closely followed by the Laplacian kernel, whereas the bilinear kernel performs significantly worse.

\section{Conclusion}\label{sec:conclusion}

In this paper, we treated the problem of identifying a nonnegative operator from given input-output trajectories. We considered a class of input-output operators defined via a feature map and an arbitrary nonnegative bounded linear operator in a ``sum of squares" fashion, such that nonnegativity of the input-output operator is guaranteed. Then, we cast the identification problem as a regularized least squares problem over a cone of nonnegative bounded linear operators. As our main result, we proved a representer theorem for this problem in the case where the input space is finite-dimensional. We also showed that the latter provides a computationally tractable solution which can be obtained by solving an associated semidefinite program. For future research,  in addition to nonnegativity, we will investigate the problem of imposing causality of the identified input-output operator.

\bibliographystyle{ieeetr}
\bibliography{../../references/kernel_refs}
\end{document}